\documentclass[12pt]{amsart}
\usepackage{epsfig}
\usepackage{amsmath}
\usepackage{amsfonts}
\usepackage{amssymb}
\usepackage{epigraph, fancyhdr}

\def\F{\mathcal F}
\def\K{\mathcal K}
\def\J{\mathcal J}
\def\L{\mathcal L} 

\def\O{\mathcal O}

\def\1{\mathbf 1}
\def\M{{\overline{\mathcal M}}}
\def\N{{\mathcal N}}

\def\QQ{\mathbb Q}
\def\ZZ{\mathbb Z}
\def\CC{\mathbb C}

\def\arm{\operatorname{\bf arm}}
\def\leg{\operatorname{\bf leg}}
\def\tail{\operatorname{\bf tail}}
\def\Res{\operatorname{Res}}

\def\det{\operatorname{det}}
\def\hat{\widehat}
\def\tilde{\widetilde}

\def\a{\alpha}

\def\t{{\mathbf t}}

\def\gs{\sigma}
\def\gl{\nu}
\def\gL{\Lambda}

\def\gS{\Sigma}

\def\lan{\langle}
\def\ran{\rangle}

\def\ft{\operatorname{ft}}

\def\td{\operatorname{td}}
\def\ch{\operatorname{ch}}

\def\fake{\operatorname{fake}}
\def\tr{\operatorname{tr}}

\renewcommand{\Delta}{\triangle}

\title[Lefschetz on $\M_{0,n}/S_n$ and adelic characterization]
      {Permutation-equivariant \\ quantum K-theory III. \\
      Lefschetz' formula on $\M_{0,n}/S_n$ \\ and adelic characterization}

\author[A. Givental]{Alexander GIVENTAL}
\thanks{This material is based upon work supported by the National 
Science Foundation under Grant DMS-1007164, and by the IBS Center for Geometry 
and Physics, POSTECH, Korea.} 

\date{June 28, 2015}

\begin{document}

\begin{abstract} 

We continue our study of the genus-$0$ permutation-equivariant quantum
K-theory of the target $X=pt$ and completely determine the ``big J-function''
of this theory. The computation is based on the application of Lefschetz' fixed
point formula to the action of $S_n$ on $\M_{0,n+1}$. It is an instance of the 
general {\em adelic characterization} (which we state at the end with reference to \cite{GiT}) of quantum K-theory for any target $X$ in terms of quantum cohomology theory. 
Yet, some simplifications of non-conceptual nature
occur in this example, making it a lucid illustration to the general theory.
\end{abstract}

\maketitle

\section*{Introduction}

By definition, the {\em big J-function} of permutation-equivariant
quantum K-theory of $X=pt$ takes an input of the form
$\t(q) =\sum_{k\in \ZZ} \t_k q^k$ with coefficients $\t_k$ from a
$\lambda$-algebra $\gL$, and assumes the value 
\[ \J_X(\t) =
1-q+\t(q)+\sum_{n\geq 2} \lan \t(L),\dots,\t(L), \frac{1}{1-qL}\ran_{0,n+1}^{S_n}.\]
In our main example, $\gL=\QQ [[N_1,N_2,\dots]]$ is the algebra of symmetric
functions, and (in fact for the formal convergence's sake) we will assume that
all $\t_k$ lie in the ideal $\gL_{+}$ of symmetric function of positive degree.
In general, we will assume that $\gL$
is equipped with an augmentation homomorphism $\gL \to \QQ$, such that
$\Psi^r(\gL_{+}) \subset (\gL_{+})^2$ for $r>1$, and that all $\t_k \in \gL_{+}$.
Moreover, we will assume that modulo any positive power of $\gL_{+}$, the
input is a Laurent polynomial in $q$. Then the value $\J_{+}(\t)$ is
a rational function of $q$ modulo any positive power of $\gL_{+}$.

We will denote by $\K$ the space of $\gL$-valued functions of $q$ rational
in this sense, and call it the {\em loop space} (writing it
sometimes as $\gL (q,q^{-1})$). It is decomposed ({\em polarized}) into
the direct sum $\K=\K_{+}\oplus \K_{-}$, where $\K_{+}$ consists of $\gL$-valued
Laurent polynomials (they can have poles only at $q=0,\infty$), and $\K_{-}$ consists of {\em reduced} rational functions regular at $q=0$ and vanishing at $q=\infty$.
We call $\K_{+}$ and $\K_{-}$ {\em positive} and {\em negative} spaces of the
polarization, and write sometimes $\K_{+}=\gL [q,q^{-1}]$ (suppressing in this
notation the aforementioned completion).

The correlators in the definition of $\J_{pt}$ lie in $\K_{-}$. Thus,
the big J-function $\t \mapsto \J_{pt}(\t)$ can be considered as a the graph of
a function $\K_{+}\to \K_{-}$ defined in the infinitesimal neighborhood
$1-q+\gL_{+} [q,q^{-1}]$ of $1-q$, called the {\em dilaton shift}.
The outcome of our computations will be the
following description of this graph as a subset in $\K$. 

\medskip

{\tt Theorem.} {\em The range of the big J-function $\J_{pt}$ in the loop
  space $\K$ is swept by the family $J_{pt}(\gl) \K_{+}, \gl\in \gL_{+}$
  of subspaces obtained from $\K_{+}$ by multiplication by the small
  J-function:
  \[ \bigcup_{\gl\in \gL_{+}}\, (1-q) e^{\text{\normalsize $\sum_{k>0} \Psi^k(\gl)/k(1-q^k)$}} \K_{+}\, .\]} 

{\tt Remark.} We'd like to warn the reader about some abuse of terminology in the expression ``is swept by the family of subspaces'' we use here and may use in the forthcoming parts of this work.  The union of subspaces is certainly greater than the the range of the J-function, e.g. because the union is a cone with the vertex at the origin, where it is singular and cannot be a range of a function. Thus, the union provides some analytic continuation of the range outside the infinitesimal $\gL_{+}$-neighborhood of $1-q$, and it would be more accurate to say that the range is {\em contained in} that union. Yet the depiction of the range as swept by varying subspaces seems more graphical, and we resort to it frequently.

\section*{Lefschetz' formula as recursion}

The theorem stated in Introduction is a consequence of Quantum
Hirzebruch-Riemann-Roch theorem improving slightly the result of
\cite{GiT} to cover the permutation-equivariant theory, and reducing in the case
of the point target space to the application of Lefschetz' fixed point formula
on $\M_{0,n+1}$.

Recall the definition of correlators in the big J-function:
\[ \J_{pt}(\t):=1-q+\t(q)+\sum_{n\geq 2}\frac{1}{n!}\tr_{h\in S_n}
H^*\left(\M_{0,n+1}; \frac{\otimes_{i=1}^n\t(L_i)}{1-qL_{n+1}}\right),\]
where the sheaf cohomology is considered as a $\ZZ_2$-graded (in the sense of linear super-algebra) $S_n$-module, with the action induced by the renumbering of the first $n$ marked points.
Recall that the trace operation acts on coefficients from $\gL$ (i.e. on tensor products of $GL_N$-modules in our main example of $\gL$) via the Adams operations:
\[ \tr_h (\gl^{\otimes n}) = \prod_r \Psi^r(\gl)^{l_r(h)},\]
where $l_r(h)$ equals the number of length-$r$ cycles in the permutation $h$.

Lefschetz' fixed point formula says
   \[ \tr_h H^*(\M_{0,n+1}; \frac{\otimes_i\t(L_i)}{1-qL})
   =  \chi \left( \M_{0,n+1}^h ; 
   \frac{\tr_h(\otimes_i\t(L_i)/(1-qL))}
        {\det_{1-h} \N^*_{\M_{o,n+1}^h}}\right),\]
where $\chi$ denotes the holomorphic Euler characteristic of the virtual bundle on the locus $\M_{0,n+1}^h$ fixed by $h$, $\tr_h (V)$ is defined on fiberwise $h$-invariant bundles by decomposition into eigen-subbundles $V_{\lambda}$ according to the eigenvalues of $h$:
   \[ \tr_h (V) := \oplus_{\lambda} \lambda V_{\lambda}, \det_{1-h}(V) =
   \tr_h \bigwedge\nolimits^{\!*} V,\]
$\N^*_{\M_{0,n+1}^h}$ is the conormal bundle to $\M_{0,n+1}^h$ in $\M_{0,n+1}$, and $L$ stands for $L_{n+1}$ (and will continue to do so till the end of this section). 

The main difficulty in the application of this formula to our situation consists in the complexity of th fixed point loci. The only way for us to
handle the problem is to apply it to the whole of $\J_{pt}$ and to all $h$, and obtain a recursion relation, which begins with examining how the symmetry $h$
acts on a genus-$0$ curve near the last, $n+1$-st marked point.
The following diagram is our book-keeping device in assessing contributions of various fixed point loci into the big J-function.

If the equivalence class in $\M_{0,n+1}$ of a stable curve is invariant under a permutation $h$ of the marked points, then $h$ induces a symmetry of the curve accomplishing the permutation, but preserving the marked point that carries the input $1/(1-qL)$. We call this marked point the {\em horn} of the $h$-invariant curve, and denote by $\zeta$ the eigenvalue of $h$ on the {\em cotangent} line to the curve at this point. 

When $\zeta=1$, we call {\em head} the maximal 
connected part of the curve containing the horn where $h$ acts as the identity. The rest of the curve decomposes into connected parts which we call {\em arms}. Thus, by definition, at the node where an arm is attached to the head, the eigenvalue $\eta$ of $h$ on the cotangent line to the arm is different from $1$. 

When $\zeta\neq 1$ is a primitive $m$th root of unity, consider the maximal 
{\em balanced} part of the curve containing the horn where the $m$-th power $h^m$ acts as the identity. A node is {\em balanced} if the eigenvalues of $h$ on the two branches are inverse. Thus, although the curve
may contain a chain of $\CC P^1$ connected at the nodes, on which $h$ acts by
multiplication by $m$th roots of $1$ (and hence $h^m$ by the identity), we take 
only the connected  part of the chain where the eigenvalues of $h$ on the 
two components connected at the nodes are equal to $1/\zeta,\zeta$ (in this
order, looking away from the horn). Thus, this chain of $\CC P^1$ (typically one $\CC P^1$) starts with the horn and ends with another fixed point of $h$, the {\em butt}, with the eigenvalue $1/\zeta$. The butt could be a non-special point, a marked point, or a node where, the {\em tail} is attached. Thus, by the definition of tail, the butt is not balanced, and hence the eigenvalue of $h$ on the irreducible component of the tail at the butt is different from $\zeta$. 

\begin{figure}[htb]
\begin{center}
\epsfig{file=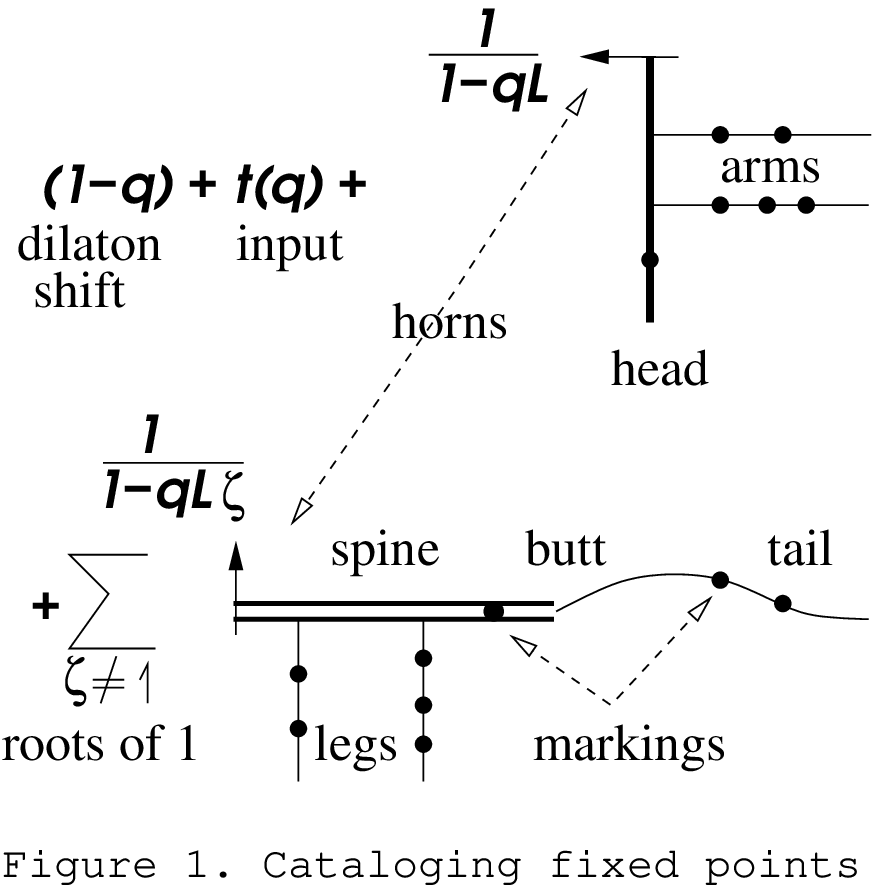}
\end{center}
\end{figure}

The (balanced chain of) $\CC P^1$ in question has a $\ZZ_m$-symmetry defined by 
$h$, and thus defines the quotient curve. We call this quotient curve the
{\em spine} (``stem'' in terminology of \cite{GiT}). In fact 
the spine comes equipped with a $\ZZ_m$-cover ramified at the the head, butt (and some nodes in the case of a chain). This means that it represents a
{\em stable map to the orbifold} $pt/\ZZ_m = B\ZZ_m$ in the sense of \cite{AGV, CR}. 

The stem can carry (unramified) marked points,
which come from symmetric configurations 
of $m$-tuples of marked points on the cover, or have nodes, which come from
$m$-tuples of symmetric nodes on the cover, where further components of the curve, cyclically permuted by $h$, but invariant under $h^m$, are attached.
In the quotient curve, these components are represented by {\em legs} attached to the spine at unramified nodes. Thus, by definition, the eigenvalue of $h^m$ on a leg at the node of attachment is different from $1$.   


Having finished our description of the diagram's elements, let us clarify that 
the diagram is meant to represent $\J_{pt}$ written as the sum of contributions 
(in the form of suitable holomorphic Euler characteristics from the right-hand-side of Lefschetz' formula) of all fixed point loci in $\M_{0,n+1}$ for all $n$, grouped by the eigenvalue $\zeta$ of the symmetry $h$ at the horn, whereas the pictures merely represent a typical appearance of the curve from the requisite group of fixed point loci. We denote by $\arm (L), \leg_{\zeta} (L), \tail_{\zeta} (L)$ the tootalities of the corresponding contributions, where $L$ represents the cotangent line bundle (over the head or spine moduli space) at the node of attachment. In fact we have
\[ \leg_{\zeta} = \tr_{h} \arm^{\otimes m} = \Psi^m(\arm),\]
since the symmetry $h$ acts by cyclically permuting the $m$ factors $\arm$ corresponding to the $m$ copies of the same leg attached to the covering 
curve.  

\smallskip

Our representation of $\J_{pt}$ has three interpretations. 

\smallskip

Firstly, one can think of it as a {\em recursion relation}. When curves in a moduli space are glued from various components, one can express the 
contribution via gluing in terms of products of contributions of the spaces of
the components. Thus, if only one manages to compute the contributions into 
Lefschetz' formula of {\em the head and spine spaces with arbitrary inputs},
this would allow one to recursively reconstruct the entire sum. Indeed, it is not an accident that the neighborhood of a horn on our diagram is shown like a curve near a node. When a typical curve in a stratum has a node, the conormal bundle to the stratum contains the deformation of smoothing the node. If $L_{+}$
and $L_{-}$ denote the cotangent line bundles at the node to the components 
of the glued curves, and $\eta_{\pm}$ the corresponding eigenvalues of the 
symmetry, the smoothing mode contributes into the denominator of Lefschetz' formula the factor $1-L_{+}\eta_{+}\eta_{-}L_{-}$. This agrees with the 
denominator $1-qL$ of the input at the horn, when $q=L_{+}\eta_{+}$ and 
$L=L_{+}\eta_{-}$. As a result, the tails and arms (and hence legs too) are themselves expressed as constituents of the total sum of $\J_X$, but simpler ones: with fewer marked points. 

Indeed, to be stable, genus-$0$ curve must carry at least $3$ special (marked or singular) points. Thus, if an arm (or tail) carries all marked point of the whole curve (a leg cannot do so, since it signifies an $m$-tuple
of identical branches with $m>1$), then no other arm (respectively, no leg) can be stable, and which in its turn shows that the head (respectively spine) will cannot be stable too.

This observation initiates a recursive procedure which in principle {\em 
completely and uniquely recovers $\J_X$ from the input $\t$ and the head and spine correlators}.

\smallskip

Secondly, one can think of our representation for $\J_X$ ``analytically'' as 
decomposition of rational functions of $q$ into elementary fractions. 
In the decomposition of $\J_{pt}$ (on the diagram) according to the eigenvalues $\zeta$, each summand represents a holomorphic Euler characteristics of a virtual bundle on a fixed point manifold. In K-theory of a manifold, the line bundle $L$ is unipotent.
The variable $q$ occurs only at the horn, and in the terms with a given $\zeta$ only in the form 
\[ \frac{1}{1-q\zeta L} = \sum_{k\geq 0} 
\frac{(q\zeta)^k (L-1)^k}{(1-q\zeta)^{k+1}}.\]
This shows that as functions of $q$ the terms with a given $\zeta$ 
are sums of simple fraction with the pole (of any order) at $q=1/\zeta$. 

\smallskip

Finally, the diagram exhibits a remarkable interpretation of localizations 
of $\J_{pt}$ at the roots of $1$ in terms of {\em ordinary} (not permutation-equivariant) quantum K-theories of $pt/\ZZ_m$, as we will explain
in the next section.

\section*{Adelic characterization}

We begin with the ordinary J-function in the quantum K-theory of the point target space: 
\[ \J^{ord}_{pt}:=1-q+t(q)+\sum_{n>1}\frac{1}{n!}\lan t(L),\dots, t(L), \frac{1}
{1-qL} \ran_{0,n+1}.\]
Note that Deligne-Mumford spaces $\M_{0,n+1}$ are manifolds (rather than orbifolds), and hence the line bundles $L_i$ as elements of $K^0(\M_{0,n+1})$ are unipotent. Consequently, the ordinary J-function takes inputs $t(q)$ from the space of power series in $q-1$, and assumes values in the spaces of Laurent series in $q-1$. We claim that:

\smallskip

\noindent {\em The Laurent series expansion of $\J_{pt}(\t)$ (with any input $\t$) near $q=1$ lies in the range of $\J_{pt}^{ord}$ (defined over the coefficient ring $\gL$).}

\smallskip

This can be seen merely by looking at the cataloging diagram and focusing on the head term.
It shows correlators of the ordinary J-function with the input
\[ t(L) = \t(L)+\arm (L), \]
where $L$ is the line bundle formed by the cotangent lines to the head at the node of attachment. It enters the arm correlators through the horn  $1/(1-LL_{-}\zeta)$ where $L_{-}$ is the cotangent line to the arm branch at the node. As we remarked earlier, the arm sums up contributions of all possible curves except those with $\zeta = 1$. But these are all the spine terms of the diagram, but only with $q$ replaced by $L$! They have no pole at $q=1$, and so can be expanded into a power series in $q-1$. Thus, we have
\[ \J_{pt}(\t)_{(1)} = 1-q+t(q)+\sum_{n>1} \frac{1}{n!} \lan t(L),\dots,t(L),
  \frac{1}{1-qL}\ran_{0,n+1} = \J^{ord}_{pt}(t),\]
where the subscript $(1)$ indicates the localization at $q=1$, i.e. expansion of rational functions of $q$ into Laurent series in $q-1$.

\smallskip

Examining the diagram again, but focusing on the terms with the horn $1/(1-q\zeta L)$ where $\zeta$ is a primitive $m$-th root of unity, we find
\begin{equation} \label{expansion_at_zeta}
  \J_{pt}(\t)_{(\zeta)}=\delta (q\zeta) +\sum_{n>0} \lan \delta (L), \Psi^m(t(L)), \dots, \Psi^m(t(L)), \frac{1}{1-q\zeta L}\ran_{0,n+2}^{spine_{\zeta}},\end{equation}
where the spine correlators are expressible somehow in terms of ordinary
quantum K-theory on $B\ZZ_m$. Here the subscript $\zeta$ means power series expansion near $q=\zeta^{-1}$.
Respectively,
\[ \delta(q) = (1-q/\zeta)+\t(q/\zeta)+\tail_{\zeta}(q/\zeta)\]
Note that the eigenvalue of symmetry $h$ on the cotangent line at the butt is $1/\zeta$, the inverse to the eigenvalue at the horn. In view of this,
the there summands in the butt input $\delta (L)$ correctly represent (in the reversed order) the possibilities that the tail is attached at the butt,
the butt being a marked point, or the butt being a regular point on the covering curve. In the last case, one is forced to declare the butt (which is a fixed point of the symmetry $h$ on the curve) a marked point on the spine. This introduces a spurious direction conormal to the fixed locus of $h$ 
(infinitesimal translation of this marked point), and the factor $1-L/\zeta$ in the numerator of Lefschetz' formula compensates for this.       

In fact equation (\ref{expansion_at_zeta}) indicates that $\J_{pt}(\t)_{(\zeta)}$
lies in a certain linear subspace in the space of Laurent series, characterized
in terms of quantum K-theory of $B\ZZ_m$. In fact the relevant stacks
of curves in $B\ZZ_m$ -- with only two ramified marked points (the horn and butt), and several (say, $l$ ) unramified --- are isomorphic (as coarse moduli spaces) to Deligne-Mumford spaces $\M_{0,l+2}$. Still, the corresponding terms
in Lefschetz' formula have complicated denominators. Namely, the $m$-fold cover of a stem curve can be deformed in $\M_{2+ml}$ away from the symmetric locus by perturbing the $m$-tuples of marked point asymmetrically, or by smoothing the nodes in an asymmetric way. We don't know an elementary way to
take this into account, but there is a machinery \cite{CGi, Co, ToT} of {\em twisted} Gromov-Witten invariants designed to handle such problems. In particular, the results \cite{ToT} by V. Tonita express the invariants of $X\times B\ZZ_m$ twisted in the required fashion in terms of the untwisted ones. The latter in their turn can be computed in terms of the ordinary K-theory of $X$ (see \cite{ACV,JK}). Below we formulate the answer for $X=pt$, and defer the general formulation to the last section. 

Let us denote by denote by $\hat{\K}$ the spaces of Laurent series in $q-1$ with coefficients in $\gL$, and by $\hat{\K}_{+}$ its subspace formed by power series, and extend the Adams operations from $\gL$ to $\hat{K}$ (and $\K$) by $\Psi^m(q):=q^m$.

\medskip

{\tt Theorem} (adelic characterization). {\em The values of the big J-function $\J_{pt}$ of the permutation-equivariant quantum K-theory of the point are completely characterized among elements of $\K$ by the requirements:

  (i) the Laurent expansion $\J_{pt}(\t)_{(1)}$ near $q=1$ is a value of $\J_{pt}^{ord}$; 

  (ii) for every primitive $m$-th root of unity $\zeta$,
  \[ \J_{pt}(\t)_{(\zeta)} (q^{1/m}/\zeta)\ \in\  \Psi^m\left(\frac{\J_{pt}(\t)_{(1)}}{1-q}\right)\, \hat{\K}_{+};\] 

  (iii) if $\zeta\neq 0,\infty$ is not a root of unity, $\J_{pt}(\t)_{(\zeta)}(q/\zeta) \in \K_{+}$.}
  
\medskip

In fact, holomorphic Euler characteristics on $\M_{0,n}$ are not too hard to compute (see \cite{YPL}), and a complete description of J-function $\J_{pt}^{ord}$ is well-known (see e.g. \cite{CGL}):
\[ \text{\em The range of $\J_{pt}^{ord}$ in $\hat{\K}$ is contained in}\ \  
  \bigcup_{\tau\in \gL_{+}} (1-q) e^{\tau/(1-q)} \hat{\K}_{+}.\]
Together with the adelic characterization theorem, this implies our Theorem from Introduction. Indeed, consider $f\in \K$ of the form
\[ f(q) = (1-q) e^{\sum_{k>0} \Psi^k(\gl)/k(1-q^k)} p(q,q^{-1}), \]
  where $p \in \K_{+}$ is a Laurent polynomial. Obviously, $f$ passes test (iii) of adelic characterization, since all poles of the exponent are at roots of unity. Next, the Laurent expansion of $f$ near $q=1$ has the form (here $\tilde{p}\in \hat{\K}_{+}$ is a power series in $q-1$):   
  \[ f(q)_{(1)} = (1-q) e^{\frac{1}{1-q}\sum_{k>0}\Psi^k(\tau)/k^2} \tilde{p}(q-1). \]
  It lies in $(1-q)e^{\tau/(1-q)}\hat{\K}_{+}$  with
    \[ \tau = \sum_{k>0} \frac{\Psi^k(\gl)}{k^2}.\]
    Thus $f$ passes test $(i)$ of adelic characterization. Finally, at a primitive $m$-th root of unity $\zeta$, we have
    \[ f_{(\zeta)} (q^{1/m}/\zeta) \in  e^{\frac{1}{1-q}
      \sum_{l>0} \Psi^{ml}(\gl)/ml^2} \hat{\K}_{+},\]
  where we used that $1/(1-q^{1/m})=m/(1-q)+\text{terms regular at $q=1$}$.
    This needs to be compared to
    \[ \Psi^m\left(\frac{f_{(1)}(q)}{1-q}\right) \hat{\K}_{+} = e^{\Psi^m(\tau)/(1-q^m)} \hat{\K}_{+}=e^{\Psi^m(\tau)/m(1-q)} \hat{\K}_{+},\]
where we used that $1/(1-q^m)=1/m(1-q)+\text{terms regular at $q=1$}$.
Since
\[ \frac{\Psi^m(\tau)}{m}=\sum_{l>0}\frac{\Psi^{ml}(\gl)}{ml^2},\]
we conclude that test (ii) has been also passed. Thus $f$ lies in the range
of the big J-function $\J_{pt}$. This range is the graph of a formal functiom
$\K_{+}\to \K_{-}$ defined in an infinitesimal neighborhood of $1-q$. Since $\Psi^1(\gl)=\gl$, it follows from the formal Implicit Fucntion Theorem that the projecion to $\K_{+}$ of the described part of the range covers the entire neighborhood, and hence the whole of the graph is thus descibed.

\section*{Generalization to arbitrary $X$}

The above formulation is a specialization to the case $X=pt$ of
the adelic characterization of genus-0 quantum K-theory of any target space $X$, which was essentially obtained in \cite{GiT} (with some aspects outsources to \cite{ToK, ToT}); the extension from the setting of that paper to the permutation-equivariant theory is immediate. 

To formulate the general result here, let us first note that for general algebraic target spaces, moduli spaces $X_{g,n,d}$ of stable maps equipped with virtual structure sheaves \cite{YPLee}, behave as ``virtual orbifolds''. Respectively, Lefschetz' formula is to be replaced by Kawasaki's Riemann-Roch formula \cite{Kaw} on $X_{g,n,d}/S_n$. The latter expresses holomorphic Euler characteristics on a given orbifold in terms of certain {\em fake} holomorphic Euler characteristics on its inertia orbifold. By fake
holomorphic Euler characteristic $\chi^{fake}({\mathcal M}; V)$ we mean $\int_{\mathcal M} \ch (V) \td (T_{\mathcal M})$, the right-hand-side of the ordinary Hirzebruch-Riemann-Roch formula. What we called above ``ordinary'' quantum K-theory of $X=pt$ in fact coincides with
the {\em fake} one, simply because Deligne-Mumford spaces $\M_{0,n}$ are manifolds. 

The big J-function in the fake quantum K-theory theory of $X$ is defined by
\[ \J_X^{fake} =1-q+t(q)+\sum_{n,d,\a}\frac{Q^d}{n!}\phi_{\a}\lan\frac{\phi^{\a}}
   {1-qL}, t(L),\dots, t(L)\ran_{0,n+1,d}^{fake},\]
   (i.e. the same way as in the ``ordinary'' quantum K-theory, but with genuine
   holomorphic Euler characteristics on $X_{g,n,d}$ replaces with the fake ones).
   Here $\{\phi_{\alpha}\}$ is a basis in $K^0(X)$, $\{\phi^{\alpha}\}$ is the dual basis with respect to the K-theoretic Poincar\'e pairing $(a,b):=\chi (X; ab)$,
   and $Q^d$ is the monomial representing the degree $d$ curves in the Novikov ring. One extends $\gL$ to include Novikov's variables, and makes it a $\lambda$-algebra by $\Psi^m(Q^d)=Q^{md}$. Thus $\J^{fake}_X$ assumes values in
   $\hat{\K}^{X}$, the space of Laurent series in $q-1$ with vector coefficients
   from $K^X:=K^0(X)\otimes \gL$, and takes inputs from $\hat{\K}^X_{+}$, the space
   of power series with such coefficients. The range of $\J_X^{fake}$ is swept
   by the family of subspaces
   \[ \bigcup_{\tau \in K^X} S_{\tau}(q)^{-1} (1-q)\K_{+}^X, \]
   where $S_{\tau}$ is a certain ``matrix'' series in $1/(q-1)$
   with coefficients which are $\gL$-endomorphisms of $K^X$ (in fact columns of $S^{-1}$ are partial derivatives of the small J-function $\J^{fake}_X(\tau)$ with respect to $\tau$). The space $S_{\tau}(q)^{-1} \K_{+}^X$ is tangent to the range of $\J_X^{fake}$ at all points of $S_{\tau}(q)^{-1} (1-q)\K_{+}^X\subset
   S_{\tau}(q)^{-1} \K_{+}^X$. With this notation, we can state:

   \medskip

   {\tt Theorem} (adelic characterization for general target $X$, cf. \cite{GiT}).
   
   {\em All values of the big J-function of permutation-equivariant quantum K-theory of $X$ are characterized by the following requirements:

     (i) $\J_X(\t)_{(1)}$ is a value of $\J_X^{fake}$;

     (ii) for every primitive $m$-th root of unity $\zeta$,
     \[ \J_X(\t)_{(\zeta)}(q^{1/m}/\zeta) \in e^{\sum_{k>0}\left(  \frac{\Psi^k(T^*_X)}{k(1-\zeta^{-k}q^{k/m})}-\frac{\Psi^{km}(T_X^*)}{k(1-q^{km})}\right)} \Psi^m \left( S_{\tau (\J_X(\t)_{(1)}}) \right) \hat{\K}^X_{+}; \]
     (iii) when $\zeta \neq 0,\infty$ is not a root of unity, $\J_X(\t)_{(\zeta)} \in \hat{\K}_{+}$.}

   \medskip
   
   By $\tau (\J_X(\t)_{(1)})$ we mean the value of the parameter of that subspace
     $S_{\tau}(q)^{-1} (1-q)\K_{+}^X$ to which $\J_X(\t)_{(1)}$ belongs. Note that the action of the Adams operation $\Psi^m$ on the operator series $S_{\tau}^{-1}$, whose ``matrix'' entries are $1/(q-1)$-series with coefficients from $\gL$, is induced not only by the $\lambda$-structure on $\gL$ (including the action on Novikov;s variables), and by $\Psi^m(q)=q^m$, but also involves conjugation by the natural automorphism $\Psi^m: K^0(X)\otimes \QQ \to K^0(X)\otimes \QQ$.

\enddocument